\documentclass[a4]{amsart}
\usepackage{amsmath,amsthm,amsfonts, amssymb, color,graphicx}
\newtheorem{theo}{\sc Theorem}
\newtheorem{lem}{\sc Lemma}
\newtheorem{rem}{\sc Remark}
\newtheorem{exa}{\sc Example}
\def\be#1{\begin{equation}\label{#1}}
\def\ee{\end{equation}}
\def\Re{{\rm Re\,}} 
\def\C{\mathbb{C}}

\def\Z{\mathbb{Z}}

\def\P{\mathfrak{P}}

\def\ds{\displaystyle}\def\ts{\textstyle}
\def\ende{{\Box}}

\def\v{\mathfrak{v}}
\def\w{\mathfrak{w}}
\def\z{\mathfrak{z}}
\def\x{\mathfrak{x}}
\def\y{\mathfrak{y}}
\def\Y{\mathfrak{Y}}
\def\barY{\tilde\Y}
\def\Ende{\hfill$\ende$}

\def\GO#1{O(|z|^{#1})}

\def\p{{\sf p}}

\def\r{{\sf r}}
\def\s{{\sf s}}
\def\t{\mathfrak{t}}
\def\tt{\mathbf{t}}\def\tt{{\sf t}}
\def\0{\mathfrak{0}}

\def\und{\quad{\rm and}\quad}
\title{\sc First order algebraic differential equations of genus zero}
\author{\sc Norbert Steinmetz}

\begin{document}
\maketitle

\bigskip

{\small\begin{center} {\sc Abstract}\end{center}
\begin{quote}  We utilise recent results about the transcendental solutions to Riccati differential equations
to provide a comprehensive description of the nature of the transcendental solutions to algebraic first order differential equations of genus zero.
\end{quote}


}

\bigskip{\it Dedicated to Professor Walter Hayman on the occasion of his 90th birthday}

\section{Introduction and Main Result}\label{INTRO}
First order differential equations
\be{FOE}P(z,w,w')=\sum_{\nu=0}^q P_{\nu}(z,w)w'^\nu=0\quad(q\ge 2)\ee
were intensively investigated by several authors \cite{BaKa1,BaKa2,AE,AE1,AAG,JM1,NSt01,KY1,KY3}.
The aim of this paper is to give a comprehensive description of the nature of the trans\-cendental solutions to equation (\ref{FOE}),
based on the re-scaling method, which was developed in the context of Riccati and Painlev\'e differential equation (\cite{NStRiccati,NStHamilton}).

\subsection{\it Equations of genus zero}
It will be assumed that $P$ is an irreducible polynomial in $\C(z)[w,w']$, hence each $P_\nu$ is a polynomial in $w$ with rational coefficients.
The hypothesis that (\ref{FOE}) has trans\-cendental meromorphic solutions is very restrictive. For example, it follows that  $\deg_w P_\nu\le 2q-2\nu$ holds;
in particular, $P_q$ is independent of $w$. Without going at length into technical and historical details -- we mention the so-called {\it Fuchsian conditions} for the absence of movable singularities other than poles, in detail deduced in Golubew~\cite{WWG} II, \S 7, see also \cite{AE}, as well as Malmquist's so-called {\it Second Theorem} \cite{JM1}(\footnote{Many authors \cite{BaKa1, BaKa2, NSt01,KY1,KY3} seemed to be unaware of Malmquist's  {\it Second Theorem}, or were in doubt about Malmquist's  reasoning. One reason might be that several elegant and transparent proofs for his {\it First Theorem} were known, but none for the second one; see also the comments in Eremenko~\cite{AE}, p.\ 62.
It is, however, not the same to ask for (necessary and sufficient) conditions for the existence of {\it one transcendental meromorphic solution} on one hand,
and for the {\it absence of movable singularities} except poles on the other.}) on the necessary and sufficient conditions that ensure the existence of transcendental meromorphic solutions
--, we will consider equations (\ref{FOE}) of {\it genus zero}. Equations of {\it genus one} will be considered in a sub-sequent paper.

For almost all parameters $z$ the algebraic curve $P(z,r,s)=0$ has some rational parametrisation $r(\tt)=\r(z,\tt),s(\tt)=\s(z,\tt)$.
By the bi-rational transformation
\be{Birational}\left.\begin{array}{rcl}w&=&\r(z,\tt)\cr w'&=&\s(z,\tt)\end{array}\right\}
\und \tt=\rho(z,w,w')\quad(\r, \s, \rho {\rm~rational~functions})\ee
equation (\ref{FOE}) is transformed into some Riccati differential equation
\be{Riccati}\tt'=a(z)+b(z)\tt+c(z)\tt^2\ee
for $\tt$, with rational coefficients $a$, $b$, and $c$. This follows from Malmquist's First Theorem applied to
$$\r_z(z,\tt)+\r_\tt(z,\tt)\tt'=\s(z,\tt).$$
Besides solutions given by (\ref{Birational}) there also occur {\it singular solutions}. They solve the {\it discriminant equation} $D(z,w)=0$ ($D$ is the discriminant of $P(z,w,w')$, regarded
as a polynomial in $w'$).

For almost every pair $(z_0,w_0)$, equation $\r(z_0,\tau)=w_0$ has $\deg\r$ distinct solutions $\tau=\tau_j$, hence from (\ref{Birational})
we obtain $\deg\r$ different solutions to equation (\ref{FOE}), defined by the initial values $w(z_0)=w_0$, $w'(z_0)=\s(z_0,\tau_j)$.
Conversely, equation~(\ref{FOE}) has $q$ solutions satisfying $w(z_0)=w_0$, $w'(z_0)=w_\nu'$, where $\omega=w_1',\ldots,w_q'$ denote the solutions to $P(z,w_0,\omega)=0$;
note that Picard's Existence and Uniqueness Theorem applies to
$$w''=\ds -\frac{P_z(z,w,w')-P_w(z,w,w')w'}{P_{w'}(z,w,w')}.$$
This shows $\deg\r=q$.

{\begin{exa}\label{Exa1}\rm The solutions to
$$zw'^2+(2zw^2-w-8z^2)w'+(zw^4-zw^2-8z^2w^2+4zw+16z^3+1)=0$$
have the form $w=\tt+z/\tt$, where $\tt$ is a solution to $\tt'=z-\tt^2$.
Like any other, this example is obtained in the following way: start with $\r=\r_1/\r_2$ and compute $\s=\s_1/\s_2=r_z+r_\tt(a+b\tt+c\tt^2)$ and the {\it resultant} $P(z,v,v_1)$ of the polynomials
$v\r_2(z,\tt)-\r_1(z,\tt)$ and $v_1\s_2(z,\tt)-\s_1(z,\tt)$ with respect to $\tt$. Then $w=\r(z,\t)$ with $\t'=a+b\tt+c\tt^2$ satisfies $P(z,w,w')=0$.\hfill$\diamondsuit$
\end{exa}}

\subsection{\it Main result}In the most simple case $q=1$, equation (\ref{FOE}) reduces to a Riccati equation (\ref{Riccati}). Based on the properties of the solutions to (\ref{Riccati}) we will prove the following
Theorem, which gives a comprehensive description of the transcendental solutions to genus-zero equations (\ref{FOE}).
For notations and results in Nevanlinna theory the reader is referred to Hayman's monograph~\cite{WKH}.

\begin{theo}\label{FOEMain}There exists some integer $n\ge -1$ and $n+2$ open sectors
$$\Sigma_{\nu}:\bar\theta_{\nu-1}<\arg z<\bar\theta_{\nu}$$ with central angle
$\frac{2\pi}{n+2}$, such that for each transcendental solution $w$ to equation (\ref{FOE}) the following holds.
\begin{itemize}
\item[{\bf a.}] $w$ has an asymptotic expansion in $\sqrt z$ (in $z$ if $n$ is even) on each $\Sigma_\nu$.
\item[{\bf b.}] Up to finitely many, the poles ($c$-points) of $w$ may be arranged in sequences $(p_k)$, each being asymptotic to one of the rays $\sigma_\nu:\arg z=\bar\theta_\nu$,
and such that $p_{k+1}=p_k\pm(\lambda\pi i+o(1))p_k^{-\frac n2}$ holds for some $\lambda\ne 0$.
\item[{\bf c.}]  $w$ has Nevanlinna characteristic $T(r,w)=Cr^{\frac n2+1}+o(r^{\frac n2})$ for some constant $C=C(w)>0$.
\item[{\bf d.}]  $w$ has at most two deficient values and rational functions; the deficiencies are integer multiples of $1/q$;
singular solutions are not deficient.\end{itemize}\end{theo}

\begin{rem}\rm The statements of Theorem can be made more explicit.
\begin{itemize}
\item[{\bf a'.}] There are two different asymptotic expansions, generically they alternate in adjacent sectors.
The sectors $\Sigma_\nu$ are called {\it Stokes sectors}.
\item[{\bf b'.}] Assuming $\deg_w P_0=2q$, generic solutions have $q$ sequences of poles $(p_k)$
that are asymptotic to each {\it Stokes ray} $\sigma_\nu:\arg z=\bar\theta_\nu$.
The sequences $(p_k)$ have counting function $\ds n(r,(p_k))=\frac{2}{(n+2)\pi|\lambda|}r^{\frac n2}+o(r^{\frac n2})$ for some $\lambda>0$.
In {\bf b.} $c$-points may be replaced by zeros of $w-\phi(z)$, where $\phi$ is any rational function.
\end{itemize}
`Generic' means `except for finitely many solutions': There are at most $n+2$ exceptional solutions $w_\mu$ such that the poles are distributed along the rays
$\arg z=\bar\theta_{\nu_h}$ with $\nu_h\in J_\mu$, ${\rm card}\;J_\mu=n+2-2d_\mu\le n$ and $\sum_\mu d_\mu=n+2$; $2d_\mu$ of the ray are `truncated'.
\begin{itemize}\item[{\bf c'.}] The integer $n$ as well as the parameter $\lambda$ and the angle $\bar\theta_0$ are only implicitly known. In Example~\ref{Exa1} we have $n=1$, $\lambda=1$,
and $\bar\theta_0=\frac{\pi}{n+2}$.
\item[{\bf d'.}] Deficient values and rational functions correspond to the asymptotic (and then convergent) series in {\bf a.}\end{itemize}
\end{rem}

The proof of Theorem~\ref{FOEMain} is based on the parametrisation (\ref{Birational}) and will be given in the next section. Of course, the parametrisation (\ref{Birational})
can be derived in particular cases only. In section~\ref{RESCALING} we will present a method how to derive the essential features--$n$ and the asymptotic series--exclusively from (\ref{FOE}).

\section{Proof of Theorem~\ref{FOEMain}}
\subsection{\it Normalisation}
Neither $\r$ and $\s$ nor the coefficients $a$, $b$, and $c(\not\equiv 0)$ in (\ref{Birational}) are uniquely determined.
For our purposes it will be convenient to choose the normal form $\tau'=P(z)-\tau^2,$
which is obtained from the original Riccati equation (\ref{Riccati}) by the simple transformation $\tt\mapsto \tau=-c(z)\tt-\frac12b(z)-\frac{c'(z)}{2c(z)}$.
The rational function
$$P=\frac14b^2-ac-\frac12b'+\frac34\Big(\frac{c'}{c}\Big)^2+\frac b2\frac{c'}c+\frac12\frac{c''}c$$
(see~\cite{HW}, p.\ 77) then satisfies
$P(z)=c_nz^n+c_{n-1}z^{n-1}+\cdots$ as $z\to\infty$, with $c_n\ne 0$ and $n\ge -1$. Although neither $P$ nor the coefficients $a,b,$ and $c$ must be polynomials,
the solutions $\tt=\rho(z,w,w')$ and $\tau$ are meromorphic in the plane.
Finally, replacing the independent variable $z$ with $\sqrt[n+2]{c_n}z$ and maintaining the notation $(z,\tt)$ then yields
\be{RiccatiNormal}\tt'=P(z)-\tt^2\quad(P(z)=z^n+a_{n-1}z^{n-1}+\cdots,~n\ge -1).\ee
By this normalisation we obtain $n$, $\bar\theta_0=\frac{\pi}{n+2}$, and $\lambda=1$.

{\begin{exa}\label{Exa1a} \rm The degrees of $a$, $b$, and $c$ in (\ref{Riccati}) may be arbitrarily large compared with $n$. This can be seen from the example (\cite{NSt01})
$$\tt'=z^{2m-1}+2z^{m-1}-(m-1)z^{m-2}+2(z^m+1)\tt+z\tt^2,$$
which may be reduced to $\tau'=1+\frac 1z-\frac3{4z^2}-\tau^2$.\hfill$\diamondsuit$
\end{exa}}

The class of equations (\ref{FOE}) is M\"oebius-invariant, that is, substitutions like
\be{Moebius}w\mapsto w_1=\frac{a(z)w+b(z)}{c(z)w+d(z)}\quad(a,b,c,d {\rm ~rational~}, ad-bc\not\equiv 0)\ee
transform equation (\ref{FOE}) into some new equation $P_1(z,w_1,w'_1)=0$ of the same type.
We will use this fact to achive
\be{P0maximal}\deg_w P_0=2q,\ee
which implies that almost all poles of $w$ are simple and $m(r,w)=O(\log r)$ holds. The Nevanlinna characteristic is $T(r,w)=N(r,w)+O(\log r)$.
The M\"obius transformations (\ref{Moebius}) also allow to transform any information about the distribution of poles of the solutions
into information about the distribution of their $c$-points and, moreover, the distribution of zeros of $w-\phi(z)$, where $\phi$ is any rational function.

The interested reader will not have any difficulty to transform the proof obtained under the
{\it special hypotheses} (\ref{RiccatiNormal}) and (\ref{P0maximal})
to the general case.

\subsection{\it Asymptotic expansions}We quote from \cite{NStRiccati} that the solutions to  (\ref{RiccatiNormal}) have asymptotic expansions
\be{Asymt}\tt(z)\sim z^{-\frac n2}\Big(\epsilon+\sum_{k=1}^\infty c_{k}(\epsilon) z^{-\frac k2}\Big)\quad(\epsilon=\epsilon_\nu(\tt)\in\{-1,1\})\ee
on the {\it Stokes sectors} $\Sigma_\nu:|\arg z-\frac{2\nu\pi}{n+2}|<\frac\pi{n+2},$
$0\le\nu\le n+1$. If $n$ is even the coefficients with $k$ odd vanish.
The coefficients $c_k$ depend on $\epsilon$, but neither on the particular solution $\tt$ nor the sector $\Sigma_\nu$.
The solution $\tt$ is uniquely determined if the asymptotic expansion holds on some open sector that contains $\overline{\Sigma}_\nu$; it then holds
on $\Sigma_{\nu-1}\cup\sigma_{\nu-1}\cup\Sigma_\nu\cup\sigma_\nu\cup\Sigma_{\nu+1}$.
Generic solutions have $\epsilon_\nu=(-1)^\nu$. To transfer the asymptotic expansions to the solutions to equation (\ref{FOE}) we need two lemmas.

\begin{lem}\label{FOlemma1}The Riccati equation (\ref{RiccatiNormal})
has at most two solutions that are algebraic at infinity. Any such solution is represented by the
(then convergent) right hand side of (\ref{Asymt}). Conversely, if the series on the right hand side converges on $|z|>r_0$, then it represents an algebraic or rational solution to (\ref{RiccatiNormal}).\end{lem}

\proof We note that asymptotic series that represent solutions on some sector are always {\it formal} solutions. If they converge at some point they converge on $|z|>r_0$ and
represent solutions that are algebraic at infinity. Conversely, every such solution $\phi$ to (\ref{RiccatiNormal}) is represented by some convergent asymptotic series (\ref{Asymt}).
Since there are only two such series, there are also at most two such solutions. \Ende

\begin{rem}\rm That there at most two such solutions also follows from the fact that the {\it cross-ratio} of any four mutually distinct solutions $\tt$, $\phi_1,\phi_2,\phi_3$ is constant. This, however, is impossible if $\tt$ is transcendental and the $\phi_\nu$ are algebraic at infinity. If there is only one $\phi$ it is rational.
In the other case, solutions $\phi_1$ and $\phi_2$ are either rational or else algebraic and analytic continuations of each other.\end{rem}

\begin{lem}\label{FOlemma2}Let $\p(z,\tau)$ be any non-constant polynomial in $\tau$ with rational coefficients, and let $\tt$ be any solution
to equation (\ref{RiccatiNormal}) with asymptotic expansion (\ref{Asymt}) on some open sector $S$. Then also $\p(z,\tt(z))$ has an asymptotic expansion in $\sqrt z$,
which is trivial (all coefficients vanish) if and only if
the equation $\p(z,\tau)=0$ has some solution $\tau=\phi(z)$ that solves the Riccati equation (\ref{RiccatiNormal}).
In this case,
\be{tminustau}\tt(z)-\phi(z)=\exp\Big(\frac{-4\epsilon }{n+2}z^{\frac n2+1}+\sum_{k=0}^{n+1} a_k z^{\frac k2}\Big)z^\kappa(1+o(1))\ee
holds as $z\to\infty$, uniformly on each closed sub-sector of $S$,  with $\kappa$ some complex constant and $\epsilon=\pm 1$ such that $\Re(\epsilon z^{\frac n2+1})>0$ on $S$.
\end{lem}

\proof It is obvious that $\p(z,\tt(z))$ has an asymptotic expansion on $S$. Vanishing of all coefficients is possible if and only if $\p(z,\tau)=0$ has an algebraic solution
$\tau=\phi(z)$ given by the series on the right hand side of (\ref{Asymt})
at $z=\infty$; by Lemma~\ref{FOlemma1}, $\phi$ solves (\ref{RiccatiNormal}). The difference $y=t-\phi$ tends to zero faster than any power $z^{-m}$ and satisfies
$$y'=-(\tt(z)+\phi(z))y=-2\epsilon_\nu z^{\frac n2}\Big(1+\sum_{k=1}^{2m} c_{k}(\epsilon_\nu) z^{-\frac k2}+O(|z|^{-m})\Big)y$$
for every integer $m$. Integrating yields
(\ref{tminustau}) on $S$, uniformly on every closed sub-sector of $S$, and $y(z)\to 0$ on $S$ requirers $\Re (\epsilon z^{\frac n2+1})>0$. \Ende

\medskip It is now easy to prove assertion {\bf a.} of Theorem~\ref{FOEMain}, and even more:

\begin{theo}\label{Expansion}Any transcendental meromorphic solution $w$ to equation (\ref{FOE}) has an asymptotic expansion on each Stokes sector $\Sigma_\nu$, except when
$\r_2(z,\tau)=0$ has a solution $\tau=\phi(z)$ given by the series on the right hand side of (\ref{Asymt}) at $z=\infty$. In this case,
$$w(z)\exp\Big(\frac{-4\ell\epsilon_\nu }{n+2}z^{\frac n2+1}+\sum_{k=0}^{n+1} \ell a_k(\epsilon_\nu) z^{\frac k2}\Big)z^{\ell\kappa(\epsilon_\nu)},$$
$\ell$ some positive integer, has an asymptotic expansion with $\Re(\epsilon_\nu z^{\frac n2+1})>0$ on $\Sigma_\nu$.\end{theo}

\proof It follows from Lemma~\ref{FOlemma2} that $\r_1(z,\tt(z))$ and $\r_2(z,\tt(z))$, hence also $\r(z,\tt(z))$, have asymptotic expansions, provided the case $\r_2(z,\tt(z))\sim 0$
(all coefficients vanish) is excluded. In this case Lemma~\ref{FOlemma2}, applied to $\p=\r_2$, gives the second statement, where $\ell$ is the multiplicity of the solution $\tau=\phi(z)$
to $\r_2(z,\tau)=0$. \Ende

\subsection{\it The distribution of poles}It follows from our hypothesis $\deg_w P_0=2q$ that $m(r,w)=O(\log r)$ and that almost all poles of $w$ are simple. They arise from
the zeros of $\r_2(z,\tt(z))$ and from the poles of $\tt$, provided $\deg_t \r_1>\deg_t\r_2$, hence $\deg_t \r_1=1+\deg_t\r_2$.
Regarding the poles of $\tt$ we recall some facts from \cite{NStRiccati}. Up to finitely many, the poles of any generic solution are arranged in $n+2$ sequences $(p_k)$
satisfying the approximate iteration scheme
$$p_{k+1}=p_k\pm(\pi i+o(1))p_k^{-\frac n2}$$
with counting function
$$n(r,(p_k))=\frac{r^{\frac n2+1}}{(n+2)\pi}+o(r^{\frac n2+1}).$$
Each such sequence, also called {\it string}, is asymptotic to some {\it Stokes ray} $\arg z=\bar\theta_\nu=\frac{2\nu+1}{n+2}\pi$.
For every pole $p\ne 0$ of $\tt$ we set $\triangle_\delta(p)=\{z:|z-p|<\delta|p|^{-\frac n2}\}$.
Then for $\delta>0$ sufficiently small the discs $\triangle_\delta(p)$ are mutually disjoint and $\tt(z)=\GO{\frac n2}$ holds outside the union $\mathcal{P}_\delta(\tt)=
\bigcup_{p\ne 0}\triangle_\delta(p)$ of these discs. We  note that the zeros $\zeta$ of $\tt$ also form strings of the same kind and are separated from the poles in the following sense:
$\liminf_{\zeta\to\infty}|\zeta|^{\frac n2}{\rm dist}(\zeta,\mathcal{P}(\tt))=\frac\pi 2$.

\begin{rem}\rm The results about the solutions to (\ref{RiccatiNormal}) are obtained with the help of the so-called re-scaling method (more in section~\ref{RESCALING}).
For any solution $\tt$ to equation (\ref{RiccatiNormal}) the re-scaled family $(\tt_h)_{|h|>1}$ of functions
$$\tt_h(\z)=h^{-\frac n2}\tt(h+h^{-\frac n2}\z)$$
is normal in the sense of Montel, and every limit function $\t=\lim_{h_k\to\infty}\tt_{h_k}$ satisfies the differential equation
$\t'=1-\t^2$ with solutions $\t\equiv\pm 1$ and $\t(\z)=\coth(\z+\z_0)$. The constant solutions give rise to the asymptotic expansions (\ref{Asymt}), while the information on the distribution of poles relies on the knowledge of the distribution of poles of the hyperbolic cotangent.\end{rem}

To prove that the zeros of $\r_2(z,\tt(z))$ are also distributed in strings we need

\begin{lem}\label{pztPoles}Let $\tt$ be any transcendental meromorphic solution to equation (\ref{RiccatiNormal}) and let $\p(z,\tau)$ be any polynomial in $\tau$ of degree  $d=\deg_\tau \p>0$
with rational coefficients, such that the solutions $\tau=\phi(z)$ to $\p(z,\tau)=0$ do not solve (\ref{RiccatiNormal}). Then
$$m\Big(r,\frac1{\p(z,\tt(z))}\Big)=O(\log r),$$
and the zeros of $\p(z,\tt(z))$ are distributed in finitely many strings, each Stokes ray attracting $d$ strings; strings of $\ell$-fold zeros will be counted $\ell$-fold.\end{lem}


\proof For $z$, $v$ and $v_1$ fixed we consider the polynomials
$Q_0(\tau)=v-\p(z,\tau)$ and $Q_1(\tau)=v_1-\p_z(z,\tau)-\p_\tau(z,\tau)(P(z)-\tau^2),$
and denote by $Q(z,s,s_1)$ the {\it resultant} of $Q_0$ and $Q_1$. Then $v(z)=\p(z,\tt(z))$ solves $Q(z,v,v')=0$ with $Q(z,0,0)\not\equiv 0$, hence $m(r,1/v)=O(\log r)$ holds by a well-known theorem due to
A.Z.\ and V.D.\ Mokhon'ko~\cite{MokMok}. To prove the second part we consider any branch $\tau=\phi(z)$ of the algebraic function $\p(z,\tau)=0$; $\phi$ is meromorphic
on the sector $S:0<\arg z<2\pi$, $|z|>r_0$ sufficiently large (note that $\arg z=0$ is not a Stokes ray). We have to discuss two cases as follows:
\begin{center}(i) $z^{\frac n2}=o(|\phi(z)|)$ and (ii) $\phi(z)=\GO{\frac n2}$ as $z\to\infty$ on $S.$\end{center}
In the first case all but finitely many of the zeros of $\tt(z)-\phi(z)$ are contained in $\mathcal{P}_\delta(\tt)$; this follows from $\tt(z)=\GO{\frac n2}$ outside $\mathcal{P}_\delta(\tt)$. We have to show that for $|p|$ sufficiently large, $\triangle_\delta(p)$
contains {\it exactly one} zero. Since $\tt$ has no zeros on $\triangle_\delta(p)$, $f(z)=\phi(z)/\tt(z)$ is regular on $\triangle_\delta(p)$ and has there exactly one zero
(namely $p$). Since $f(z)\to\infty$ as $p\to\infty$, uniformly on $\partial\triangle_\delta(p)$, Rouch\'e's Theorem applies to $f$ and $f-1$, hence $f-1$ and $\tt-\phi$
have exactly one zero on $\triangle_\delta(p)$.

In the second case we re-scale along any sequence $(\zeta_k)$ of zeros of $\tt-\phi$ to obtain the initial value problem $\t'=1-\t^2$, $\t(0)=\lim_{\zeta_k\to\infty}\zeta_k^{-\frac n2}\phi(\zeta_k)$
for the limit function $\t=\lim_{\zeta_k\to\infty}\tt_{\zeta_k}$. If $\t(0)\ne\pm 1$ it follows that the zeros of $\tt-\phi$ form strings, again with $\zeta_{k+1}=\zeta_k\pm(\pi i+o(1))\zeta_k^{-\frac n2}$, and exactly one in each Stokes direction.

If, however, $\t(0)=\pm 1$, hence $\phi(z)=\pm z^{\frac n2}+\cdots$, but $\varrho=\phi'-P+\phi^2\not\equiv 0$, holds we will consider $u=1/(\tt-\phi(z))$ and the corresponding differential equation $u'=1+2\phi(z)u+\varrho(z)u^2,$
which may be transformed into normal form
$$v'=P^*(z)-v^2,$$
now with $P^*(z)=P(z)-2\phi'(z)+\frac34\big(\frac{\varrho'(z)}{\varrho(z)}\big)^2+\phi(z)\frac{\varrho'(z)}{\varrho(z)}-\frac12\frac{\varrho''(z)}{\varrho(z)}=z^{n}+\cdots$. This proves that the zeros of $\tt-\phi$, which coincide with the poles of $v$, also form strings of the same kind.
\Ende

\medskip We have thus proved assertion {\bf b.} of Theorem~\ref{FOEMain}.

{\begin{exa}\label{Exa3a}\rm Let $\tt$ be any generic solution to $\tt'=z^2-\t^2$. Then
$$\ds w=\frac{t^4}{(t-z)(t-2z)(t-z^2)}$$
solves some equation (\ref{FOE}) of degree four. It has four different types of poles distributed in strings along the Stokes rays $\arg z=(2\nu+1)\frac{\pi}4$ and corresponding to
\begin{enumerate}
\item  the poles $p$ of $\tt$; they form the set $\mathcal{P}(\tt)$.
\item the zeros $\zeta$ of $\tt(z)-z^2$; they are contained in $\mathcal{P}_\delta(\tt)$, exactly one belongs to $\triangle_\delta(p)$ for $|p|$ sufficiently large; actually
$\zeta=p+o(|p|^{-1})$.
\item the zeros $\tilde \zeta$ of $\tt(z)-2z$; re-scaling along any sequence $(\tilde \zeta_k)$ leads to the initial value problem $\t'=1-\t^2$, $\t(0)=2$, hence
$\t(\z)=\coth(\z+\frac12\log 3)$. The pole of $\t$ closest to $\z=0$ is $-\frac12\log 3$, hence the pole of $\tt$ closest
to $\tilde \zeta_k$ is $p_k=\tilde \zeta_k-(\frac12\log 3 +o(1))\tilde \zeta_k^{-1}$ by Hurwitz' Theorem.
The poles $\tilde p$ form strings that are `parallel' to the strings of the first kind.
\item the zeros $\hat \zeta$ of $\tt(z)-z$; re-scaling along any sequence $\hat \zeta_k)$ leads to the initial value problem $\t'=1-\t^2$, $\t(0)=1$, hence
$\t(\z)=1$. The poles of $\tt$ are `invisible' from $\hat \zeta_k$ in the metric $ds=|z||dz|$, that is, $|\hat \zeta_k|{\rm dist}(\hat \zeta_k,\mathcal{P}(\tt))$ tends to
infinity as $k\to\infty$.~\hfill$\diamondsuit$\end{enumerate}
\end{exa}}

\subsection{\it The Nevanlinna characteristic}Since $w$ has $q(n+2)$ strings of poles, the total number of poles on $|z|<r$ is
$$n(r,w)=\ds q(n+2)\frac{r^{\frac n2+1}}{(\frac n2+1)\pi}+o(r^{\frac n2+1})=\frac 2\pi qr^{\frac n2+1}+o(r^{\frac n2+1}),$$
and we obtain
$\ds T(r,w)=N(r,w)+O(\log r)=\frac{4q}{(n+2)\pi}r^{\frac n2+1}+o(r^{\frac n2+1}).$
This proves Theorem~\ref{FOEMain} {\bf c.} for generic solutions. In the exceptional cases there is an additional factor $1-\frac{2d(\tt)}{n+2}.$

\subsection{\it Deficient values and rational functions}Let $\psi$ be any rational function or constant. The already mentioned  A.Z.\ and V.D.\ Mokhon'ko-Theorem~\cite{MokMok}
yields
$$m\Big(r,\frac1{w-\psi}\Big)=O(\log r)\und \delta(\psi,w)=\liminf_{r\to\infty}\frac{m\big(r,\frac1{w-\psi}\big)}{T(r,w)}=0,$$
provided $P(z,\psi(z),\psi'(z))\not\equiv 0$. On the other hand, if $\psi$ solves (\ref{FOE}) but is not singular, the algebraic equation $\r(z,\tau)-\psi(z)=0$,
equivalently $\p(z,\tau)=\r_1(z,\tau)-\phi(z)\r_2(z,\tau)=0$
has solutions $\tau=\phi(z)$ that also solve equation (\ref{RiccatiNormal}). By Lemma~\ref{FOlemma1}, any such $\phi$ is given by the (now convergent) series on the right hand side of (\ref{Asymt}).
This, in particular, implies that there are at most two deficient rational functions or constants of this kind.
Now $\p$ factors into $\p_1\p_2$, such that all solutions to $\p_1(z,\tau)=0$, but none to $\p_2(z,\tau)=0$
also solve the Riccati equation (\ref{RiccatiNormal}). From
$$w-\psi(z)=\frac{\r_2(z,\tt(z))}{\p_2(z,\tt(z))}\p_1(z,\tt(z))$$
and
$$\ds m\Big(r,\frac{\r_2(z,\tt(z))}{\p_2(z,\tt(z))}\Big)+m\Big(r,\frac{\p_2(z,\tt(z))}{\r_2(z,\tt(z))}\Big)=O(\log r)$$
(which follows from Lemma~\ref{pztPoles}, since none of the solutions to equation $\r_2(z,\tau)=0$ also solves (\ref{RiccatiNormal})) it then follows that
$$m\Big(r,\frac1{w-\psi}\Big)=m\Big(r,\frac1{\p_1(z,\tt(z))}\Big)+O(\log r).$$
The Uniqueness Theorem for the initial value problem $\tt'=P(z)-\tt^2$, $\tt(z_0)=\phi(z_0)$, where $\phi$ is any solution to $\p_1(z,\tau)=0$ then shows that
$\p_1(z,\tt(z))$ has only finitely many zeros, hence we obtain
$$\begin{array}{rcl}
\ds m\Big(r,\frac1{w-\psi}\Big)&=&\ds T(r,\p_1(z,\tt(z))+O(\log r)\cr
&=&\ds\deg_\tau\p_1 T(r,\t)+O(\log r)=\frac{4\deg_\tau\p_1}{(n+2)\pi}r^{\frac n2+1}+o(r^{\frac n2+1})\end{array}$$
and $\delta(\psi,w)=\frac{\deg_\tau\p_1}{\deg_\tau\r}=\frac{\deg_\tau \p_1}{q}$. Singular solutions $\psi$, however, have $\r_1\equiv 1$ 
and $\delta(\psi,w)=0$. \Ende

{\begin{exa}\label{Exa4}{\bf (\cite{NSt01})} \rm Equation
$(w'-2b(z)w)^2=4w(a(z)+c(z)w)^2$ $(a,b,c$ polynomials, $ac\not\equiv 0)$ arises from the Riccati equation
$\tt'=a(z)+b(z)\tt+c(z)\tt^2$
by the simple transformation $w=\t^2$. In any case, the values $0$ and $\infty$ are completely ramified for $w$.
For $a+b+c\equiv a+2b+4c\equiv 0$, say, $w$ has two deficient values: $\delta(1,w)=\delta(4,w)=\frac12$.\hfill$\diamondsuit$
\end{exa}}

\section{Re-scaling}\label{RESCALING}
\subsection{\it Algebraic differential equations}Suppose $w$ is any meromorphic solution to some algebraic differential equation
\be{ADE}Q(z,w,w',\ldots,w^{(n)})=0.\ee
To obtain the (essential) properties of $w$, set
\be{wsubh}w_h(\z)=h^{-\alpha}w(h+h^{-\beta}\z),\ee
$w'_h(\z)=h^{-\alpha-\beta}w(h+h^{-\beta}\z)$ etc
to obtain $Q(\z,h^{\alpha}w_h,h^{\alpha+\beta}w_h',\ldots,h^{\alpha+n\beta}w_h^{(n)})=0.$
Taking the limit
$\mathfrak{Q}(\x_0,\x_1,\ldots,\x_n)=\lim_{h\to\infty}h^{-m}Q(h,h^\alpha\x_0,h^{\alpha+\beta}\x_1,\ldots,h^{\alpha+n\beta}\x_n)$
for suitably chosen $m$ then yields the autonomous equation
\be{ADErescaled}\mathfrak{Q}(\w,\w',\ldots,\w^{(n)})=0\ee
for $\w=\lim_{h\to\infty}w_{h}$. Apart from the fact that the real parameters $\alpha$ and $\beta$ are arbitrary, the method is by no means justified.
Nevertheless it can be justified if the functions $w_h$ form a normal family in the sense of Montel.

\subsection{\it A normality criterion}Normality of this family may be characterised by the growth of some generalisation to {\it spherical derivative}
as follows.

\begin{lem}{\rm(\cite{NStPainleve})} Normality of {\it any} re-scaled family $(w_h)_{|h|>1}$ defined by (\ref{wsubh}) is equivalent to
\be{NH}\limsup_{z\to\infty}\frac{|w'(z)||z|^{\alpha-\beta}}{|z|^{2\alpha}+|w(z)|^2}<\infty\ee
\end{lem}

\begin{rem}\rm The re-scaling method was introduced in \cite{NStPII, NStRiccati, NStHamilton} in the context of various analytic differential equations. It was inspired by the well-known Zalcman Re-scaling Lemma~\cite{LZ1,LZ2} and  Yosida's work \cite{Yosida2}. Let $\alpha$ and $\beta>-1$ be real parameters. The class $\barY_{\alpha,\beta}$ consists of all meromorphic functions $f$ such that the family $(w_h)_{|h|>1}$ of functions (\ref{wsubh})
is normal on $\C$ in the sense of {Montel}, and all limit functions  $\w=\lim_{h_k\to\infty}w_{h_k}$ are $\not\equiv\infty$, at least one of them being non-constant.
If, in addition, {\it all} limit functions are non-constant, then $w$ is said to belong to the {\it {Yosida} class}
$\Y_{\alpha,\beta}$. The class $\Y_{0,0}$ was introduced by {Yosida}~\cite{Yosida2} (denoted $A_0$ there), and for arbitrary real parameters by the author~\cite{NStYosida};
it is universal in the sense that it contains all limit functions $\w=\lim_{h_n\to\infty}w_{h_n}$
for $w\in\Y_{\alpha,\beta}$. Instead of $w$ meromorphic in the plane one could also consider $w$ meromorphic on some sector $S$
(and restrict the sequences $(h_n)$ to arbitrary closed sub-sectors of $S$).\end{rem}

{\begin{exa}\label{Exa3}\rm The solutions to (\ref{RiccatiNormal}) belong to the class $\barY_{\frac n2,\frac n2}$, the (components of the) solutions to the Hamiltonian system $p'=-q^2-zp-a$, $q'=p^2+zq+b$
belong to $\barY_{1,1}$, and the first, second, and fourth Painlev\'e transcendents belong to the classes
$\barY_{\frac 12,\frac 14}$, $\barY_{\frac 12,\frac 12}$, and $\barY_{1,1}$, respectively (see \cite{NStRiccati, NStHamilton,NStPainleve}).\hfill$\diamondsuit$\end{exa}}

\subsection{\it Application to first order differential equations} In the present case of (\ref{FOE}) the formal re-scaling process yields so-called Briot-Bouquet equations
\be{FOErescaled}\mathfrak{P}(\w,\w')=0.\ee
The solutions to (\ref{FOErescaled}) belong to the class $W$ (like Weierstrass, notation introduced by Eremenko), which consists or rational, trigonometric, and elliptic functions. In our case (genus zero) elliptic functions do not occur. Given any equation (\ref{FOE}) of genus zero, neither the parametrisation (\ref{Birational}) nor the Riccati equation  (\ref{RiccatiNormal}) are at hand.
Thus the problem arises how to determine the parameters $n$ and $\theta_0$ as well as the asymptotic expansion immediately from (\ref{FOE}),
and also possible values of $\alpha$ and $\beta=\frac n2$, if any, such that $w\in\barY_{\alpha,\beta}$.\medskip

{\it a. Asymptotic expansions.} To determine the potential leading term of asymptotic expansion $w\sim az^{\frac m2}+\cdots$ exclusively from equation (\ref{FOE}) (it follows from (\ref{Birational})
that $m$ must be some integer), consider $P(x,y,{\ts\frac m2} yx^{-1})=\sum_{\nu=0}^{2q}(A_\nu+o(1))x^{k_\nu}y^\nu$
and apply the Newton-Puiseux method to the simplified equation
\be{ReducedEquation}\sum_{\nu=0}^{2q}A_\nu x^{k_\nu}y^\nu=0.\ee
As $x\to\infty$, the solutions have leading terms $a_jx^{\rho_j}$ ($a_j\ne 0$), and the potential leading terms of the asymptotic expansions are among the terms $a_jz^{\rho_j}$
with $2\rho_j=m_j\in\Z$. Some of these terms may also belong to singular solutions.\medskip

{\it b. The parameters $\alpha$ and $\beta$.} To determine the possible values of $\alpha$ and $\beta$ we will just consider equations (\ref{FOE}) such that (\ref{FOErescaled}), which is obtained by a formal limiting process,
has maximal degree $\deg_{\w'}\P=q$. This is a reasonable postulate since we want to deduce all relevant properties from (\ref{FOErescaled}). We assume $P_q(z,w)\equiv 1$.
For $\x$ and $\y$ fixed, consider
$$\Phi(h,\x,\y)=h^{-q(\alpha+\beta)}P(h,h^\alpha\x,h^{\alpha+\beta}\y)=\mathfrak{P}(\x,\y)+\phi(h,\x,\y)$$
with $\P(\x,\y)=\y^q+\cdots$ and $\deg_\y\phi<q$. Then $\alpha$ and $\beta$ can be adjusted in such a way that $\phi(h,\x,\y)$ tends to zero as $h\to\infty$.
Of course, the procedure is not unique, and we aim to choose $\beta$ as small as possible (in order that the `local unit discs' $|z-p|<|p|^{-\beta}$ are as large as possible).\medskip

{\it c. Proof of $w\in\barY_{\alpha,\beta}$.} Having determined the possible parameters one has to prove $|w'|=O\big(|z|^{\beta-\alpha}(|z|^{2\alpha}+|w|^2)\big)$. This may be done by using well known estimates for the roots of an ordinary equation $x^q+p_{q-1}x^{q-1}+\cdots +p_0$ applied to (\ref{FOE}), where $P$ has to be regarded  as a polynomial in $w'$.

\subsection{\it Examples} We will now give some examples to illustrate the method. Non-trivial examples necessarily look quite complicated.
{\begin{exa}\label{ExaWeberHermite}\rm  Consider $z^2w'^2+P_1(z,w)w'+P_0(z,w)=0$ with
\be{BSPWH}
\begin{array}{rcl}P_1(z,w)&=&(2z-2z^3)w-\frac1{4}(2-z^2)w^2, {\rm~and}\cr
P_0(z,w)&=&2z^5w+(1+\frac{31}4z^4)w^2+(4z+5z^3)w^3-\frac14(2+3z^2)w^4.\end{array}\ee
a. The reduced equation~(\ref{ReducedEquation}) is given by $yx^2(32x^3+124x^2y+80xy^2-12y^3)=0$ with solutions $0, 8x, -\frac13x, -x$. The first pair corresponds to the singular solutions $w=0$ and $w=8z$ (the discriminant of $P$ is  $D(z,w)=16z^4w(w-8z)((2+7z^2)w+4z^3)^2$), while the second pair determines the principal terms of the asymptotic expansions $w\sim-\frac13z+\cdots$ and $w\sim-z+\cdots$\medskip

b. For any choice of $\alpha$ and $\beta$, the principal part of $\Phi(h,\x,\y)$ has the form
$$\ts\y^2-(2h^{1-\beta}\x-\frac14h^{\alpha-\beta}\x^2)\y+2h^{3-\alpha-2\beta}\x+\frac{31}4h^{2-2\beta}\x^2+5h^{1+\alpha-2\beta}\x^3-\frac34h^{2\alpha-2\beta}\x^4$$
Obviously, $\beta\ge 1$ is necessary. Choosing $\beta=1$, the terms $-\frac14h^{\alpha-1}\x^2\y$ and $2h^{1-\alpha}\x$ enforce $\alpha=1$
and
$\P(\w,\w')=\ts\w'^2-(2\w-\frac14\w^2)\w'+2\w+\frac{31}4\w^2+5\w^3-\frac34\w^4=0,$
with solutions $\w=0, 8, -\frac13,-1$, and $\ds\w=\frac{\coth^2\z}{\coth\z-2}$.\medskip

c. To prove $|w'|=O(|z|+|w|^2)$ we use the well-known upper bounds for the zeros of ordinary polynomials to obtain
$$\begin{array}{rcl}
|w'|&=&O(\max\{|P_1(z,w)/z^2|,|P_0(z,w)/z^2|^{\frac 12}\})\cr
&=&O(\max\{|z||w|, |w|^2|, |z|^{\frac 32}|w|^{\frac12}, |z||w|^{\frac 32}\}).\end{array}$$
Applying various H\"older inequalities then shows $|w'|=O(|z|^2+|w|^2)$.
Re-scaling along any sequence of poles yields non-constant limit functions.
To prove that $w_{h_k}\to \infty$ as $h_k\to\infty$ is impossible we consider the  differential equation for $\v=1/\w$,
$\ts\v'^2+(2\v-\frac14)\v'+2\v^3+\frac{31}4\v^2+5\v-\frac34=0,$
which has no trivial solution. This proves  $w\in\barY_{1,1}$.\medskip

{\it Distribution of poles.} Almost all poles of $w$ are simple and are distributed in strings asymptotic to the rays $\arg z=(2\nu+1)\frac\pi2$. The poles of $w$ occur in pairs $p$ and $\tilde p=p+\frac12\log 3+o(1))p^{-1}$ (this following from $\coth\frac12\log 3=2$.\hfill$\diamondsuit$
\end{exa}}

{\begin{exa}\label{ExaWeberHermite_b}\rm The same procedure applied to
$$w'^2-z^2(4w-w^2)w'+4z^4w+(8z^2-z^2)w^2+(4-6z^2)w^3-(1-z^2)w^4$$
yields $\alpha=0$ and $\beta=2$, which, however, doesn't reflect the properties of the transcendental solutions
$w=\ds\frac{\tt^2}{\tt-1}$ with $\tt'=z^2-\tt^2$ and $n=2$, hence $\beta=1$.\hfill$\diamondsuit$\end{exa}}

{\begin{exa}\label{Exa7}\rm Equation $z^2w'^2+[(2z-4z^2)w-(1-z)w^2]w'+$
$$[4z^3w+(1-z^2+8z^3)w^2+(4z-6z^2+4z^3)w^3-(1-z+z^2)w^4]=0$$
is obtained from  $\tt'=z-\tt^2$ by the transformation
$\ds w=\frac{\tt^2}{\tt-z}.$
From our method we obtain $\beta=\frac12$ (this is not surprising since $\tt\in\barY_{\frac12,\frac12}$) and $\alpha=0$. Again from the differential equation it follows that $|w'|=O(|z|^{\frac12}(1+|w|^2)$, hence the family of functions
$w_h(\z)=w(h+h^{-\frac12}\z)$ is normal ($\alpha=0$ and $\beta=\frac 12$). The limit equation is $\w'^2=-4\w(\w+1)^2$ with solutions $\w=0, -1$ and $\w=-\coth^2\z$. Like in Example the limit $\w=\infty$  does not occur, and-constant limit functions are obtained by re-scaling along any
sequence $h_k$ such that $w(h_k)=c\ne 0,-1$, say. This proves $w\in\barY_{0,\frac12}$.
Since non-constant solutions $\w$ have double poles, any pole $p$ of $w$ is accompanied by a pole $\tilde p=p+\epsilon(p) p^{-\frac12}$
with $\epsilon(p)\to 0$ as $p\to\infty$.
More precisely, if $\tilde p$ has residue $1$ and $p$ has residue $-\frac{p^2-p+1}{p^2+p}$ (obtained from the differential equation), re-scaling along $p=p_k$ gives
$$w_{p_k}(\z)=\frac{\epsilon(p_k)p_k^2-\epsilon(p_k)(p_k-1)+(2p_k-1)\z}{\sqrt{p_k}(p_k+1)\z(\z-\epsilon(p_k)}+{\rm ~bounded~function~on~} |\z|<\delta.$$
In order that $w_{p_k}$ tends to $\w=-\coth^2\z=1/{\z^2}+O(1)$ it is necessary and sufficient that $\epsilon(p_k)\sim -p_k^{-\frac 12}$ and $\tilde p_k\sim p_k-p_k^{-1}$.
Leaving the singular solutions $w=4z$ and $w=0$ aside (discriminant $D(z,w)=z^4w(w-4z)((1-z+2z^2)w+2z^2)^2$), the constant solution $\w=-1$ leads to two different asymptotic expansions
$w\sim -1\pm z^{-\frac 12}-z^{-1}\pm\cdots$,
which may be computed immediately from the differential equation. \hfill$\diamondsuit$\end{exa}}

\bigskip
\noindent {\small\it Norbert Steinmetz.  Institut f\"ur Mathematik.
Technische Universit\"at Dortmund. \\ D-44221 Dortmund, Germany.
E-mail: stein@math.tu-dortmund.de.\\
Web: http://www.mathematik.tu-dortmund.de/steinmetz/}

\end{document}